**Enabling Dynamic Process Simulators to Perform Alternative Tasks: A Time-stepper Based Toolkit for Computer-Aided Analysis**

C.I.Siettos[1], C. C. Pantelides[2] and I. G. Kevrekidis[1,3]

[1]Department of Chemical Engineering, Princeton University, Princeton, 08544, NJ, USA
[2]Centre for Process Systems Engineering,
Department of Chemical Engineering and Chemical Technology,
Imperial College of Science, Technology and Medicine, London, SW72BY, UK
[3]also PACM and Department of Mathematics, Princeton University

**Abstract**

We discuss computational superstructures that, using repeated, appropriately initialized short calls, enable temporal process simulators to perform alternative tasks such as fixed point computation, stability analysis and projective integration. We illustrate these concepts through the acceleration of a gPROMS-based Rapid Pressure Swing Adsorption simulation, and discuss their scope and possible extensions.

## Introduction

Good process design is both a challenge and a necessity in the chemical process industry, and scientific computation plays a vital role in this endeavor. Computational engineers make increasing use of advances in numerical analysis or algorithm development that have the potential to improve current modeling practices. The purpose of this note is to bring to the attention of computational design practitioners a set of "numerical superstructures" that enable state-of-the-art *process temporal simulation codes* to perform a number of tasks (such as fixed point location, continuation, stability/bifurcation analysis and possibly even controller design and optimisation) that they may have not been explicitly designed for. These superstructures, which one might group under the label "Numerical Analysis of Legacy Codes" (e.g. [1,2,3]) have their conceptual roots in large-scale, matrix-free iterative linear algebra / subspace iteration methods. We will illustrate this time-stepper based approach to computer-aided process analysis by enabling a particular commercial modeling tool (gPROMS, [4]) to compute, and analyze the stability of, cyclic steady states (CSS) of a particular process (RPSA) faster than through direct simulation. We will comment on the programming issues involved in wrapping the appropriate computational superstructure around gPROMS. We will conclude with a discussion of potentially useful variants of this "computational enabling" theme, including what we term "coarse projective integration".

**Direct Temporal Simulation vs. Other Computational Tasks**

A (good !) temporal simulation code evolves a model of the system on the computer emulating the way the system evolves in nature: operating parameters are given, initial conditions are prescribed and the evolution of the system state in time is recorded. Good temporal simulation codes often embody many man-years of effort, and – even though they contain the best knowledge of the process available - can be difficult to maintain and modify, especially as the original programmers (whether in industry, national laboratories or research groups) move on.

More specifically, consider that the physical model comes in the form of a well-posed continuum partial differential equation set $U_t = \Phi(U;\lambda)$ for the system state $U(x,t)$ depending on parameter(s) $\lambda$. For certain design tasks (e.g. the location of steady states), temporal simulation (direct integration) is an acceptable numerical option: if the steady state in question is (globally) stable, repeated calls to the integration routine over, say, a fixed time interval *T* will eventually lead the trajectory to its neighborhood. This use of the integrator can be thought of as a successive substitution iterative scheme of the form $U^{(n+1)} = \Phi(U^{(n)};\lambda)$ where $\Phi(U;\lambda)$ represents the result of the integration over a period *T* with initial condition *U.*

Alternative algorithms, however, such as Newton-type iterations, augmented by pseudo-arclength continuation [6], are better suited –given good initial guesses- to compute steady states and follow their dependence on parameters. Augmenting the "natural" right hand side of a dynamic problem with criticality or optimality conditions gives rise to algorithms that accurately locate instability boundaries or local extrema. These new, augmented sets of equations are based on the physical model *and on the type of task we want to perform*. A single Newton solution of the system $\{U - \Phi\ (U;\lambda) = 0;\ \det[I - \Phi_U(U;\ \lambda)] = 0\ \}$ would, for example, accurately locate an ignition point $(U^*, \lambda)$ of the dynamical model of a combustion reactor much more economically than direct simulation. The latter would require extensive integrations over several initial conditions and at several parameter values in order to accurately bracket the critical parameter value. The same limitation would, of course, hold for physical experiments.

In summary, the problem we are facing is as follows: in many industrially relevant situations we already have good temporal simulators ("timesteppers", as we will refer to them from now on). However, we need information from the model that is not easily obtained through temporal simulation. We do not wish to write new code from scratch in order to obtain this information; we would rather exploit the existing simulator by enabling it, through a computational superstructure, to perform tasks it was not originally designed for [1,2,3,6]. This is the realm of "numerical analysis of legacy simulators": the construction of algorithms that have a two-tier structure. At the inner level, the algorithm calls the timestepper as a black box subroutine for relatively short integration times, and with appropriately chosen initial conditions. The results of these short calls are used to *estimate* "on demand" [5] quantities that the outer level code needs: residuals, the action of local Jacobians and local Hessians etc. Finally the outer level code uses this information to perform the task we want, such as iteration of a contraction mapping for locating steady states, the subsequent design of a stabilizing controller, or an optimisation step.

In the original work of Shroff and Keller [6], from which our inspiration originates, the inner iteration was a short-term integration step performed by an existing timestepper. It is important to notice, however, that this two-level approach is also applicable to cases where the inner simulator performs a complex set of tasks (e.g. multistage dynamic integration involving PDAEs over an entire operation cycle, with intermediate discrete decisions). The acceleration and computer aided analysis techniques we describe below may, under appropriate conditions, apply equally successfully. A case in point is the use of such time-stepper based techniques for the so-called coarse integration/bifurcation analysis of systems described by microscopic/stochastic simulators [1,2,3,7,8,9].

This two-level approach is analogous to the framework of large-scale iterative linear algebra [10]. There, the inner iteration is a simple matrix-vector product computation while the outer algorithm solves linear equations, or eigenproblems (Krylov-type methods, Arnoldi methods etc.). In a matrix-free context, nearby function evaluations help *estimate* (as opposed to directly evaluate) the necessary matrix-vector products [11].

**The gPROMS modeling software.**

gPROMS (Process Systems Enterprise, [4]) is a state-of-the-art commercial modeling environment. gPROMS models comprise descriptions of the transient behavior of the underlying physical system expressed in terms of mixed systems of ordinary differential and algebraic equations (DAEs) or integro-partial differential algebraic equations (IPDAEs). The equations can also incorporate intrinsic discontinuities in the physical behavior (e.g. those associated with the appearance or disappearance of thermodynamic phases or changes in flow regimes) described in terms of general State-Transition Networks. Complex external manipulations and disturbances (e.g. operating procedures for start-up, shut-down, emergency handling etc.) that affect the system can also be modeled in detail.

gPROMS is a multipurpose tool allowing diverse activities to be applied to the same system model. These include steady-state and dynamic simulation and optimisation, parameter estimation from steady-state and dynamic experiments, and the optimal design of such experiments. A combination of advanced symbolic and numerical solution techniques is provided for this purpose.

One of the main advantages of gPROMS in the context of the present work is its open software architecture. This allows the gPROMS engine (gSERVER, Process Systems Enterprise, [12]) to be embedded within an external code. In particular, it is possible for any gPROMS dynamic simulation of arbitrary complexity to be invoked repeatedly as a procedure, with its initial conditions and, possibly, other parameters and characteristics changing from one invocation to another.

**A representative example: Rapid Pressure Swing Adsorption, RPSA.**

Periodic adsorption processes (e.g. pressure or temperature swing adsorption) play a key role in industrial gas separation in the iron and steel, refinery, chemical and petrochemical industries. Optimizing their cyclic steady state (CSS) performance has important economic implications, and fast computation of the cyclic steady state is a vital part of this undertaking [13,14,15]. A frequent characteristic of dynamic simulators of PSA is, however, the excruciatingly slow convergence to the final CSS, which in turn becomes a bottleneck for the design procedure. We will illustrate the use of a two-tier algorithm (the Recursive Projection Method of Shroff and Keller) built around gPROMS for efficient CSS location and stability analysis of RPSA.

RPSA operation involves a cycle with two steps: pressurization by feed gas, when the impurities are adsorbed, and counter-current depressurization with internal purging. The process under study concerns the production of oxygen-enriched air from a nitrogen and methane mixture in a packed bed of zeolite 5A. During the first step of a cycle, pressurized air is fed at the bottom of the bed. Nitrogen is selectively adsorbed while the oxygen-enriched product is drawn from the top of the bed. At the second step, the feed is interrupted and the bed is depressurized simultaneously from both ends. The adsorbed nitrogen is released and leaves the bed from the bottom while oxygen-enriched product continues to be obtained from the top.

The mathematical model of the adsorption bed involves mass balances for the gas and solid phases. It is assumed that all adsorption parameters (i.e bed fraction, bed bulk density and particle size) are constant over the bed which operates isothermally. The governing equations and their boundary conditions [15] constitute a mixed set of nonlinear partial differential and algebraic equations (PDAE) that require efficient numerical techniques for their solution.

The columns was spatially discretized using a centered finite difference method of order two and a total of 90 discretization intervals, while the absolute and relative integration tolerances were set to $10^{-6}$. The evolution of the system through its first 4 cycles of operation is shown in figure 1a.

After an operation sequence of many periodic cycles, the system approaches its CSS, at which the conditions in the bed at the start and end of each cycle are identical. If $U^{(n)}$ denotes the state of the system at the start of the first step of cycle *n*, then the evaluation of the final state of the system involves the integration of the model equations over a single cycle, and establishes the transformation $U^{(n)} \rightarrow U^{(n+1)} \equiv \Phi(U^{(n)}; \lambda)$. The approach of both the model and the actual process to the CSS is sluggish, as illustrated in Fig. 1b. In fact, the fixed point $U^f$ is reached ($\left\| U^{(n)} - U^{(n+1)} \right\| \approx 10^{-5}$) after ~ 4000 operating cycles. Linearized stability of the fixed point $U^f$ is described by the eigenvalues (Floquet multipliers) of the linearized mapping $X^{(n+1)} = \Phi_U X^{(n)}$ where $X \equiv U - U^f$, and $\Phi_U$ is the monodromy matrix.

Since direct simulation converges slowly, we implemented the Recursive Projection Method (RPM) of Schroff and Keller around the gPROMS timestepper to accelerate it. The method uses the timestepper to perform three distinct tasks: (a) approximate in an iterative manner, through judicious choice of initial conditions, the slow stable (or slightly unstable) eigenspace $P$ of $\Phi_U$, which is assumed to be of a small dimension; (b) eliminate (through integration) the fast decaying modes of the solution, corresponding to fast, strongly stable eigenvalues, in $Q$, the orthogonal complement of $P$; and (c) in the adaptively identified low-dimensional subspace $P$, accelerate convergence to the fixed point through an approximate Newton method (see Fig. 2). The procedure enables the integrator to find even mildly unstable steady states to which it would never normally converge. It is straightforward to incorporate pseudo-arclength continuation in the RPM algorithm to automatically allow tracing of the solution branch through singular (turning) points (see for example the original work of Keller, or 21). Since no turning points exist in the vicinity of our nominal computation, we did not perform such continuations here.

Like all Newton-type schemes, the RPM procedure requires an initial guess that is sufficiently close to the solution. In this case, we obtain this by an initial integration of the system over 800 cycles. With this initial guess, RPM converges to the final CSS (within an error of $10^{-6}$) using 40 individual cycle integrations, each initialized at an appropriate initial state $U$. The Newton steps are performed in the resulting slow subspace, whose dimension was only 3 (compared to the size of the full model which involves 366 differential and 459 algebraic variables).

The apparent acceleration achieved by RPM over the direct simulation approach is a factor of 80 (40/3200). Even including the initial 800 cycles required to obtain a good initial guess, the acceleration is still a factor of 5 or 6. Due to the slow dynamics, successful initialization of the RPM also requires here a good initial estimate of the slow subspace. This "good initialization" overhead will be incurred much less frequently in a continuation/optimization context.

Figure 3a compares the RPM-converged CSS concentration profiles for both gas and solid phases with the results of long time integration. Since upon convergence of RPM we have an estimate of the "slow Jacobian" of the process, we can use it to estimate the slow Floquet multipliers (shown in figure 3b). As expected by the slow convergence, these critical modes are inside *but very close* to the unit circle. Upon convergence, an iterative Arnoldi algorithm [16] exploiting the same timestepper was used to validate these slow Floquet multiplier estimates. The inset in Fig. 3b shows four of the leading Floquet multipliers computed by the Arnoldi procedure. Note that our particular Arnoldi procedure, pioneered in [16], for stability determination using subspace identification with full nonlinear timesteppers also constitutes a timestepper "enabling technology".

This example goes beyond the original Shroff and Keller RPM, which was designed for differential equations and steady state solutions, in two ways. What is located here is the fixed point of the stroboscopic map of a periodically forced system – in effect, we are finding a limit cycle through shooting [17]. Furthermore, this example involved not simply differential, but *differential algebraic* timesteppers. The RPM- and Arnoldi-based linear algebra operations were performed in the subspace of the *differential variables*, $U$. Given the values for these variables at the start of a cycle, gPROMS automatically computes consistent values of the algebraic variables by applying a Newton-type procedure to the model's algebraic equations [18]. Thus, the gPROMS environment allows us to implement our software superstructure working transparently on the differential variables alone.

In the work of Keller and his students von Sosen and Love [19,20], systems of differential-algebraic equations were first explicitly converted to systems of ordinary differential equations, on which RPM was applied. A brief discussion of the issues arising in the use of RPM-type methods for differential-algebraic systems of equations can be found in [21] for a particular reactor modeling example.

**Discussion**

The goal of enabling existing legacy timesteppers is a very practical one. Instead of writing "the best possible code" for a task from scratch, and spending time to debug and validate it, we can sometimes quickly and efficiently use a validated form of our model in the form of an existing timestepper. Although the two-tiered codes will almost certainly, and possibly dramatically, underperform codes written expressly for a given task, we save in "startup programming costs" and in time-to-results.

The superstructure is common across different problems – the only modification necessary is to allow the timestepper to be called as a black box, input-output subroutine. The gPROMS open architecture makes the implementation of timestepper-based algorithms relatively straightforward. The inner iteration is simply a gPROMS-based dynamic simulation executed with initial conditions and, possibly values of other parameters, provided by the outer algorithm (RPM, Arnoldi etc.). The latter can be implemented either in an external high-level language (e.g. FORTRAN) program that runs the gPROMS engine (gSERVER) as an embedded application or directly within the gPROMS TASK language itself. Both approaches have been successfully pursued in this work.

As has already been mentioned, our RPM procedure was applied in the subspace of the differential variables, $U$ for DAE systems of index-1. An alternative approach would be to apply it in the full space $(U, V)$, where $V$ denotes the algebraic variables, in conjunction with any systematic method for obtaining a set of consistent initial values that satisfy all algebraic constraints $f(U,V) = 0$. In fact, the approach applied in the context of the RPSA example is equivalent to determining $V$ such that $f(U^*,V) = 0$ where $U^*$ are the values determined by the RPM procedure at each outer iteration. However, one could alternatively determine values $U$ and $V$ by projecting the values $U^*$ and $V^*$ suggested by the RPM procedure onto the manifold $f(U,V) = 0$, i.e. by solving the optimization problem min $(\|U-U^*\| + \|V-V^*\|)$ subject to $f(U,V)=0$. A large number of other alternatives is likely to exist in non-trivial problems (e.g. [22,23]); for example, one could fix a subset of the differential and the algebraic variables at the values suggested by the RPM procedure, and then determine consistent values of all remaining variables. In general, all of these alternatives will result in the correct operation of the outer algorithm. However, they may result in widely differing computational cost of the inner iteration overhead as they affect the ease of performing the consistent initialization step of the DAE integration.

A remarkable byproduct of subspace-based procedures like RPM (and more generally, Newton-Picard type algorithms, [24,25]) is the slow Jacobian of the inner iteration, the timestepper. For systems with separation of time scales, this is precisely the information required to perform, for example, stabilizing controller design [26,27]. Comparably, one can use the timestepper to identify the action not just of slow Jacobians, but also of slow Hessians, and thus create two-tier timestepper-based optimization algorithms.

Beyond the "legacy code" justification, there are several other motivations for the construction of two-tier timestepper based algorithms. The timestepper can be arbitrarily complex, involving, for example, discrete manipulations and discontinuities, or indeed multiscale models. In the latter case, as demonstrated originally in [3] (see also [1] and references therein), the RPM procedure operates on a set of macroscopic variables, while the timestepper evolves a microscopic/stochastic description of the system. In this context, the consistent initialization of the microscopic system involves a process of "lifting" (or "disaggregation") of the values of the macroscopic variables into one or more consistent microscopic realizations. At the end of its computation, the timestepper returns macroscopic variable values obtained via a "restriction" (or "aggregation") procedure applied to the final values of the microscopic variables.

Moreover, the "inner iteration" does not have to be a temporal integration; it can be a contraction mapping for finding a steady state (e.g. inexact Newton based on GMRES, [28]) or it can be a local (e.g. conjugate gradient type) optimization step. In general, techniques like RPM can conceivably accelerate/stabilize any iterative procedure given the conditions articulated in [6], and integration is but one such example.

We conclude with a small additional example of how such two-tier methods can be used to assist the computational analysis of engineering problems for which timestepper subroutines are available. An inspection of the transient evolution of the PSA problem in Fig. 1b clearly shows two time scales corresponding to fast periodic oscillation on the one hand, and slow evolution of the "average" state on the other. Slow, averaged equations have traditionally been derived for periodically forced oscillators in nonlinear mechanics [29,30]. Although no such equation exists for the RPSA case, direct time integration provides an approximate timestepper for this unavailable equation. More specifically, if we consider the state at the end of each cycle as representative of this slow envelope equation, we can accelerate its evolution as follows: Starting with a given initial condition we integrate for a number of cycles. We then use the system state $U$ at the ends of the last few cycles to *estimate* the right-hand-side of the equation for the slow evolution. Given this estimate, we can perform a long, explicit (e.g. Euler) time step *for this envelope equation*. This idea, for general problems with a separation of time scales, lies behind the "projective integration" methods of Gear and Kevrekidis [2,31,32]**.** A quick implementation of this idea in gPROMS results in the remarkable savings shown in Fig. 1b; this is of course a problem ideally suited for the approach, due to the large gap between the "slow" and "fast oscillatory" modes of the envelope equation.

The equation for which the projective integration step is implemented is the (explicitly unavailable) envelope equation for the rapidly oscillating problem. In principle, one could also work with the averaged equation (the "slow" equation for the average state over one forcing period). The method can also be applied to rapidly oscillating *autonomous* systems; but now the averaging –or the envelope- would have to be taken between successive crossings of a Poincaré map. Of course, the consistent initialization issues discussed earlier also arise in this projective integration context, both for the "envelope" and for the "averaged" implementation. A particular case would be Hamiltonian systems, for which one would have to impose conservation of energy (or additional integrals, if applicable) as a correction to the "unavailable averaged" or "unavailable envelope" equation projective predictor. The numerical analysis of such projective integration methods, for which a separation of time scales argument is necessary, is an ongoing research subject in the literature (see for example [31,32,33]).

In conclusion, one may find better things to do with an existing validated process timestepper than just integrate: added value can be extracted through computational superstructures motivated by subspace iteration methods from large-scale iterative linear algebra. These methods, in the presence of sufficient time scale (and concomitant space scale) separation in the original problem may significantly benefit its computer-assisted analysis. It is possible that the separation of time scales arises not in the equation itself, but in appropriate coarse-grained views of the original problem, like averaged or "envelope" equations demonstrated above, or like the coarse molecular dynamics, coarse Brownian Dynamics and coarse Kinetic Monte Carlo approaches in [1,7,8,9].

We should not give the impression that these methods are all new – although considerable original work has been performed recently along these lines, from the Shroff and Keller RPM to our projective and telescopic projective integrators and the "coarse" analysis of microscopic simulators. For example, the successful application of Broyden's method for the CSS computation of periodically operated processes [34,35,36,37,38,39] falls under this umbrella. Scientists in several fields are becoming increasingly conscious of "Newton-Picard" type methods. As far back as 1991, Sotirchos [40] accelerated periodic chemical vapor infiltration computations by, in effect, exploiting the separation of time scales of an "unavailable envelope equation". The purpose of this note is to make design practitioners aware of a potentially useful shortcut with low programming overhead, enabling them to perform diverse calculations by exploiting existing and validated legacy dynamic codes. These methods may not, in general, be computationally optimal for the task in question; but they are very flexible, general purpose and versatile. As such, this numerical enabling technology may catalyze the computer-assisted analysis of many complex processes for which validated timesteppers exist.

**Acknowledgements**

This work, originally presented at the 2001 AIChE meeting in Reno (paper no. [41]), was partially supported through AFOSR (Dynamics and Control, Drs. Jacobs and King) and UTRC (IGK, CIS) and United Kingdom's Engineering and Physical Sciences Research Council (EPSRC) under Platform Grant GR/N08636 (CCP). The effort of enabling gPROMS timesteppers through FPI involved extensive collaboration with Drs. N. Bozinis (at Imperial College) and K. Theodoropoulos (now at UMIST), which we acknowledge. We would also like to acknowledge extensive discussions and collaboration over the years with D. Roose and K. Lust, of Leuven, on their Newton-Picard large limit cycle computation methods. While this paper was being written, we learned of the independent work on the application of Newton-Picard computations to PSA processes by van Noorden, Veduyn Lunel and Bliek [42], also reported in the PhD thesis of T. L. van Noorden [43] this last June.

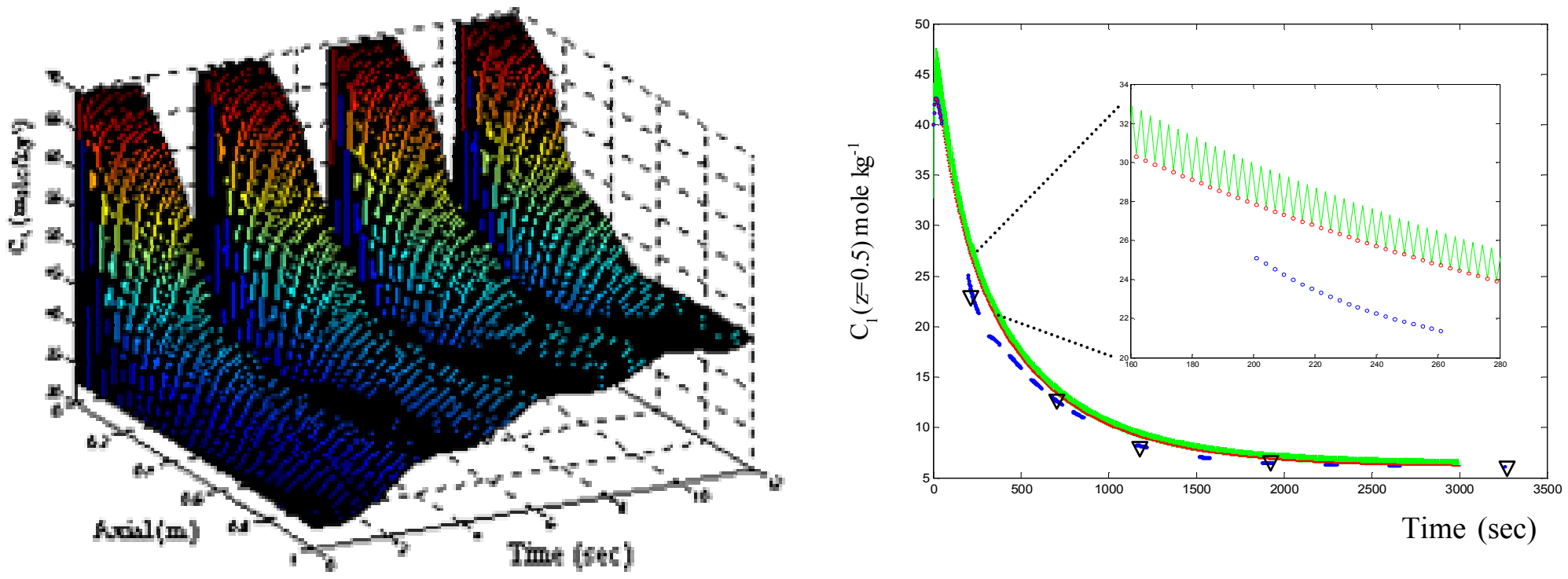

**Figure 1.** (a) Gas phase concentration for the first 4 cycles of operation. (b) Gas phase concentration at the middle of the reactor ($\nabla$: "coarse" integration).

**(I)**

*Choose* $\varepsilon$, $n_{max}$

**Give Initial conditions:** $U^{(0)}$, $\lambda$

*Do While* ($n<n_{max}$)

$U^{(n+1)} = \Phi(U^{(n)};\lambda)$ gPROMS

*If* $\left\| U^{(n+1)} - \Phi(U^{(n)};\lambda) \right\| < \varepsilon$

$U^{(final)} = U^{(n+1)}$

*stop*

*Endif*

*End Do*

**(II)**

**For n = 1 until Convergence Do**

1. Solution Decomposition $U^{(n+1)} = p^{(n+1)} + q^{(n+1)}$;

   $p = f(p,q;\lambda) \equiv P\Phi(p+q;\lambda) \quad q = g(p,q;\lambda) \equiv Q\Phi(p+q;\lambda)$

2. Newton iteration: $\tilde{p}^{(n+1)} = p^{(n)} - (f_p(p^{(n)},q^{(n)};\lambda) - \mathrm{I})^{-1}(f(p^{(n)},q^{(n)};\lambda) - p^{(n)})$

3. Picard iteration: $\tilde{q}^{(n+1)} = g(p^{(n)},q^{(n)};\lambda)$ gPROMS

**End For**

$U^{(final)} = p^{(final)} + q^{(final)}$

**Figure 2**. RPM as a pseudo code.

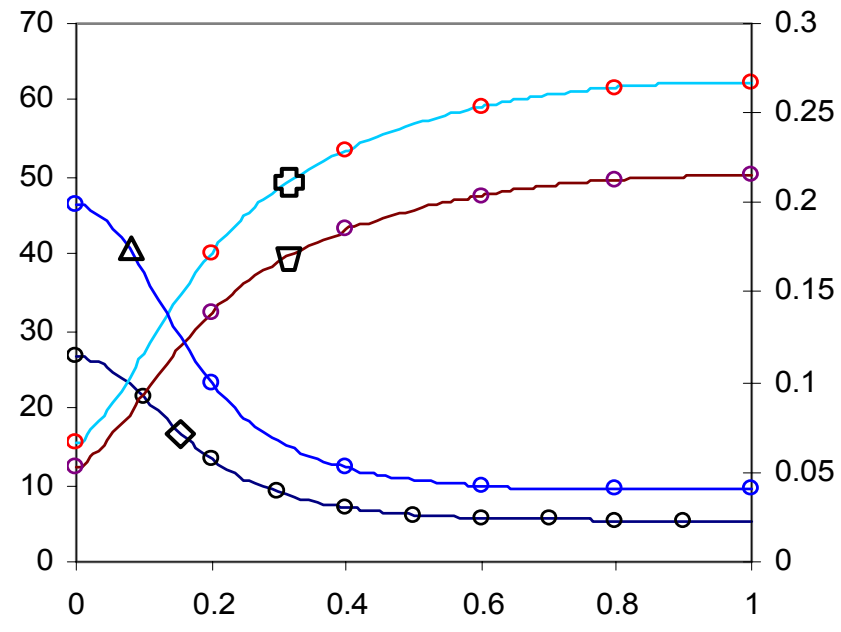

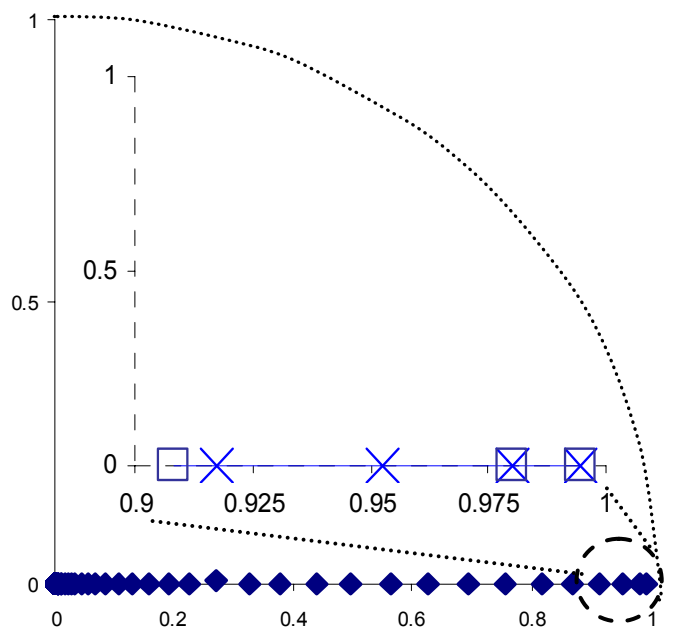

**(a)** **(b)**

**Figure 3.** (a) Solid phase (△, ▽) and gas phase (✣, ◇) concentration profiles as computed with RPM (solid line), profiles obtained through time integration after 4000 cycles of operation (O). (b) Computed Floquet multipliers: RPM (□) , Arnoldi (X ).